\documentclass[a4paper,12pt]{amsart}
\usepackage{mathrsfs}
\usepackage{amsfonts}
\usepackage{amssymb}
\usepackage{ifthen}
\usepackage{graphicx}
\usepackage[usenames]{color}
\nonstopmode \numberwithin{equation}{section}
\newtheorem{thm}{Theorem}
\newtheorem{cor}{Corollary}
\newtheorem{lem}{Lemma}

\newtheorem{claim}{Claim}
\newtheorem{conj}[equation]{Conjecture}

\theoremstyle{definition}
\newtheorem{defn}{Definition}
\newtheorem{case}{Case}
\newtheorem{subcase}{Subcase}

\newtheorem{examp}{Example}
\newtheorem{prob}[equation]{Problem}
\newtheorem{ques}[equation]{Question}
\newtheorem{rem}{Remark}

\newcounter {own}
\def\theown {\thesection       .\arabic{own}}

\newenvironment{pf}[1][]{%
 \vskip 3mm
 \noindent
 \ifthenelse{\equal{#1}{}}%
  {{\slshape Proof. }}%
  {{\slshape #1.} }%
 }%
{\qed\bigskip}

\newcounter{alphabet}
\newcounter{tmp}

\newcommand{\IN}{{\mathbb N}}
\newcommand{\IC}{{\mathbb C}}
\newcommand{\ID}{{\mathbb D}}

\newcommand{\Aut}{{\operatorname{Aut}}}

\def\be{\begin{equation}}
\def\ee{\end{equation}}

\newcommand{\bee}{\begin{enumerate}}
\newcommand{\eee}{\end{enumerate}}

\newcommand{\blem}{\begin{lem}}
\newcommand{\elem}{\end{lem}}
\newcommand{\bthm}{\begin{thm}}
\newcommand{\ethm}{\end{thm}}
\newcommand{\bcor}{\begin{cor}}
\newcommand{\ecor}{\end{cor}}
\newcommand{\beg}{\begin{examp}}
\newcommand{\eeg}{\end{examp}}
\newcommand{\begs}{\begin{examples}}
\newcommand{\eegs}{\end{examples}}
\newcommand{\bdefe}{\begin{defn}}
\newcommand{\edefe}{\end{defn}}
\newcommand{\bprob}{\begin{prob}}
\newcommand{\eprob}{\end{prob}}
\newcommand{\bques}{\begin{ques}}
\newcommand{\eques}{\end{ques}}
\newcommand{\bei}{\begin{itemize}}
\newcommand{\eei}{\end{itemize}}
\newcommand{\bca}{\begin{case}}
\newcommand{\eca}{\end{case}}
\newcommand{\bsca}{\begin{subcase}}
\newcommand{\esca}{\end{subcase}}

\newcommand{\bcl}{\begin{claim}}
\newcommand{\ecl}{\end{claim}}

\newcommand{\bcon}{\begin{conj}}
\newcommand{\econ}{\end{conj}}
\newcommand{\bcons}{\begin{conjs}}
\newcommand{\econs}{\end{conjs}}
\newcommand{\bprop}{\begin{propo}}
\newcommand{\eprop}{\end{propo}}
\newcommand{\br}{\begin{rem}}
\newcommand{\er}{\end{rem}}
\newcommand{\brs}{\begin{rems}}
\newcommand{\ers}{\end{rems}}
\newcommand{\bo}{\begin{obser}}
\newcommand{\eo}{\end{obser}}
\newcommand{\bos}{\begin{obsers}}
\newcommand{\eos}{\end{obsers}}
\newcommand{\bpf}{\begin{pf}}
\newcommand{\epf}{\end{pf}}
\newcommand{\ba}{\begin{array}}
\newcommand{\ea}{\end{array}}
\newcommand{\beq}{\begin{eqnarray}}
\newcommand{\beqq}{\begin{eqnarray*}}
\newcommand{\eeq}{\end{eqnarray}}
\newcommand{\eeqq}{\end{eqnarray*}}

\newcommand{\ds}{\displaystyle}

\begin{document}
\bibliographystyle{amsplain}
\title {Properties of Normal Harmonic Mappings}

\author[H. Deng]{Hua Deng}
\address{Hua Deng, Department of Mathematics, Hebei University,
Baoding, Hebei 071002, People's Republic of China.}
\email{1120087434@qq.com}

\author[S. Ponnusamy]{Saminathan Ponnusamy}

\address{Saminathan Ponnusamy, Department of Mathematics, Indian Institute of
Technology Madras, Chennai-600 036, India. }
\email{samy@iitm.ac.in}

\author[J. Qiao]{Jinjing Qiao ${}^{~\mathbf{*}}$}
\address{Jinjing Qiao, Department of Mathematics, Hebei University,
Baoding, Hebei 071002, People's Republic of China.}
\email{mathqiao@126.com}

\subjclass[2000]{Primary:  30D45, 31A05; Secondary: 30G30, 30H05}
\keywords{ Normal functions, Normal harmonic mappings, Spherical derivative, Maximum principle.
\\
${}^{\mathbf{*}}$ Corresponding author}


\begin{abstract}
In this paper, we present several necessary and sufficient conditions for a harmonic mapping to be normal.
Also, we discuss maximum principle and five-point theorem for normal harmonic mappings.
Furthermore, we investigate the convergence of sequences for sense-preserving normal harmonic mappings and
show that the asymptotic values and angular limits are identical for normal harmonic mappings.
\end{abstract}

\thanks{The research of this paper is supported by NSF of Hebei Science Foundation (No. A2018201033).}

\maketitle
\pagestyle{myheadings} \markboth{H. Deng, S. Ponnusamy and J. Qiao}{Properties of normal harmonic mappings}

\section{Introduction and Main results}\label{sec-1}
Let $\mathbb{D}= \{z\in \IC:\, |z|<1\}$ denote the unit disk in the complex plane $\mathbb{C}$.
A function $f$ meromorphic in $\mathbb{D}$ is called a $normal\:function$
if the family $\mathfrak{F}=\{f\circ\varphi:\, \varphi\in \Aut (\mathbb{D})\}$ is a normal family,
where $\Aut (\mathbb{D})$  denotes the class of conformal automorphisms of $\mathbb{D}$ (cf. \cite{LV}).
Normal functions were first studied by  Yosida \cite{Yo}. Subsequently, Noshiro \cite{No} gave a characterization of
normal functions by showing that a meromorphic function $f$ is normal if and only if
\be\label{eq1.1}
\sup_{z\in\ID}(1-|z|^2)f^{\#}(z)<\infty,
\ee
where $f^{\#}$ denotes the spherical derivative of $f$ given by $f^{\#}(z)={|f'(z)|}/{(1+|f(z)|^2)}.$
The condition \eqref{eq1.1} is equivalent to say that $f$ is Lipschitz when regarded as a function from
the hyperbolic disk $\mathbb{D}$ into the extended complex plane endowed with the chordal distance
(cf. \cite{LV}) which is defined as follows: The  $chordal\;distance$ $\chi(a,b)$ between the complex values $a$ and $b$, considered as points
on the Riemann sphere, is given by
\be\label{DPQ2-e q1}
\chi(a,b)=\left\{
\begin{array}{rcl}
0                                                    &       &\mbox{if $a=b$,}\\
\ds \frac{|a-b|}{\sqrt{1+|a|^2}\sqrt{1+|b|^2}}           &       &\mbox{if $a\neq\infty\neq b$,}\\
\ds \frac{1}{\sqrt{1+|a|^2}}                             &       &\mbox {if $a\neq\infty=b$.}
\end{array} \right.
\ee
Normal functions play important roles in studying properties of meormorphic functions, specially the behaviour in
the boundary of meormorphic functions. Many results have appeared in the literature, see, for example,
\cite{La1, La2, LV, LP, Po, Ya}.
The main focus in this article is to extend a number of results
from theory of  analytic functions to the case of planar harmonic mappings.

Let  $\Omega$ be a simply connected domain in $\IC$. A harmonic mapping $f$ on $\Omega$ is a complex-valued function which
has the canonical decomposition $f=h+\overline{g}$, where $h$ and $g$ are analytic in $\Omega $
and $g(z_0)=0$ at some prescribed point $z_0\in \Omega$. We recall that (see \cite{Le}) a necessary and sufficient condition for a
complex-valued harmonic mapping $f=h+\overline{g}$ is locally univalent and
sense-preserving in $\ID$ is that red $h'(z)\neq 0$ and the Jacobian $J_f(z)$ is positive in $\ID$, where $J_f(z) =|h'(z)|^2-|g'(z)|^2.$


A harmonic mapping $f=h+\overline{g}$ in $\mathbb{D}$ is red said to be normal if
$$\sup_{z_1\neq z_2}\frac{\chi(f(z_1), f(z_2))}{\rho(z_1, z_2)}<\infty,
$$
where $\rho(z_1, z_2)$ denotes the hyperbolic distance  between two points $z_1$ and $z_2$ in $\ID$, that is,
$$\rho(z_1, z_2)=\frac{1}{2}\log \left ( \frac{1+r}{1-r}\right ),\quad r=\left |\frac{z_1-z_2}{1-\bar{z_1}z_2}\right |.
$$

Following the idea of Colonna \cite{Co}  on harmonic Bloch functions,  Arbel\'{a}ez et al. \cite{AHS} studied normal
harmonic mappings and established some necessary conditions for a harmonic mapping to be normal. We begin
with the following equivalent definition (see \cite[Proposition 1]{AHS}).

\bdefe\label{de1.1}
A harmonic mapping $f=h+\overline{g}$ in $\mathbb{D}$ is  said to be normal  if
$$\sup_{z\in\mathbb{D}}(1-|z|^2)f^{\#}(z)<\infty,
$$
where
$$ f^{\#}(z)=\frac{|h'(z)|+|g'(z)|}{1+|f(z)|^2}.
$$
\edefe


Following the investigation of \cite{AHS},   we continue in this paper the study of normal harmonic mappings.
First we extend the theorem of Lohwater-Pommerenke \cite[Theorem 1]{LP}  to the case of normal harmonic mappings, in the
following form.

\bthm\label{thm1}
A non-constant function $f$ harmonic in $\mathbb{D}$ is normal if and only
if there do not exist sequences $\{z_n\}$ and $\{\rho_n\}$ with $z_n\in \mathbb{D}$, $\rho_n>0$,
$\rho_n\rightarrow 0$ as $n\rightarrow \infty$, such that
\be\label{eq2.0}
\lim_{n\to \infty}f(z_n+\rho_n\zeta)=F(\zeta)
\ee
locally uniformly in $\mathbb{C}$, where $F$ is a non-constant harmonic mapping.\ethm

It is worth pointing out that the important use of \cite[Theorem 1]{LP} was to prove the five-point theorem
due to Lappan \cite[Theorem 1]{La1} which asserts that a function $f$  meromorphic
in $\ID$ is normal if $\sup_{z\in f^{-1}(E)}f^{\#}(z)(1-|z|^2)$ is bounded for some five-point set
$E\subset f(\mathbb{D})$. Being stated Theorem \ref{thm1}, it is natural to ask whether this result continues
to hold in the case of harmonic mappings to be normal.

\bthm\label{thm2}
Let $E$ be any set consisting of five complex numbers, finite or infinite. If $f$ is a sense-preserving harmonic mapping in $\mathbb{D}$ such that
$$ \sup_{z\in f^{-1}(E)}(1-|z|^{2})f^{\#}(z)<\infty,
$$
then $f$ is a normal harmonic mapping.
\ethm

Our next result is a natural generalization of \cite[Lemma 1]{La2} from the case of normal meromorphic functions to the case of harmonic mappings.

\bthm\label{thm3}
Let $K$ be a positive real number and let $f$ be a normal harmonic mapping in $\mathbb{D}$.  Then for
each positive integer $n$, there exists a constant $E_{n}(f,K)$  satisfying the inequality
$$ (1-|z|^{2})^n(|h^{(n)}(z)|+|g^{(n)}(z)|)\leq E_{n}(f,K)
$$
for each $z\in \mathbb{D}$ and that $|f(z)|\leq K$.
\ethm

Theorem \ref{thm3} actually characterizes sense-preserving normal harmonic mappings. For if $f$ is
a sense-preserving harmonic mapping in $\mathbb{D}$ which is not normal, then by Theorem \ref{thm2} for each fixed
$K>0$ and for each value $w$ such that $|w|<K$, we have
$$\sup_{z\in f^{-1}(w)}(1-|z|^{2})f^{\#}(z)=\infty
$$
with at most four exceptions for $w$.

Moreover in \cite{La0}, Lappan showed that a meromorphic function $f$ is normal if and only if
$$\lim_{n\to \infty}f(z_n)=\lim_{n\to \infty}f(w_n)
$$
for all sequences $\{z_n\}$ and  $\{w_n\}$ in $\mathbb{D}$ such that $\rho(z_n, w_n)\to 0$ as $n\to \infty$.
This result has a natural analog for normal harmonic mappings.

\bthm\label{thm4}
Let $f=h+\overline{g}$ be a harmonic mapping in $\mathbb{D}$ such that either $h$ or $g$ is bounded.
Then  $f$ is normal in $\mathbb{D}$ if and only if for each pair of sequences $\{z_n\}$ and  $\{w_n\}$
of $\mathbb{D}$ such that $\rho(z_n, w_n)\to 0$ as $n\to \infty$,
the convergence of $\{f(z_n)\}$ implies the convergence of $\{f(w_n)\}$ with the same limit.
\ethm

Theorem \ref{thm4} characterizes normal harmonic mappings which is indeed a generalization of \cite[Theorem 1.3]{AMR}
for normal functions.

\bthm\label{thm5}
Let $f$ be harmonic in $\mathbb{D}$ and $0<p<\infty$. Then $f$ is normal if and only if
$$\sup_{z,w\in \mathbb{D}, z\neq w}\frac{\chi(f(z), f(w))}{|z-w|}|1-\overline{w}z|^{1-\frac{2}{p}}(1-|w|^2)^{\frac{1}{p}}(1-|z|^2)^{\frac{1}{p}}<\infty.
$$
\ethm

The proof of this result is similar to the proof of \cite[Theorem 1.3]{AMR} and so, we omit its proof. Note that
the case $p=2$ of Theorem \ref{thm5} gives a compact and useful form for a harmonic function to be normal.
%

The maximum principle for normal functions is established in \cite{LV} (see also \cite[Theorem 9.1]{Po1}), as a
generalization of the classical maximum principle for analytic functions since there
is no assumption on $|f(z)|$ with $z$ belonging to some subarc of the boundary. We next consider the maximum
principle for normal harmonic mappings, and get a harmonic analog of \cite[Theorem 9.1]{Po1}, which is indeed a
generalization of the classical maximum principle for harmonic mappings.

\bthm\label{thm4.1}
Let $f=h+\overline{g}$ be harmonic  in $\mathbb{D}$ and
\be\label{eq4.0}
\sup_{z\in\mathbb{D}}(1-|z|^2)\frac{|h'(z)|+|g'(z)|}{1+|f(z)|^2}\leq\alpha<\infty.
\ee
Let $G$ be a domain with $\overline{G}\subset\mathbb{D}$ that lies in the lens-shaped domain of angle
$\beta~(0<\beta<\pi)$ cut off from $\mathbb{D}$ by the circular arc $B$ (see Figure \ref{maximum}). We suppose that
\be\label{eq4.1}
|f(z)|\leq\delta< \delta_0
\ee
for $z\in\partial G\backslash B$, where
$\delta_0=\frac{1}{\kappa}\left (1+\sqrt{1+\kappa^2}\right )\exp \left [-\sqrt{1+\kappa^2}\right ]$
with $\kappa=\frac{\alpha\beta}{\sin\beta}$. Then
\be\label{eq4.2}
|f(z)|\leq\eta ~\mbox{ for $z\in G$,}
\ee
where $\eta=\eta(\delta,\alpha,\beta)$ is the smallest positive solution of
\be\label{eq4.3}
\delta=b(\eta), \quad  b(t)=t\exp\left (-\frac{\kappa}{2}\Big (t+\frac{1}{t}\Big )\right ).
\ee
\ethm

It is a simple exercise to see that the function $b(t)=t\exp\left (-\frac{\kappa}{2}\left (t+\frac{1}{t}\right )\right )$
is increasing for $0<t<t_0$ and decreasing for $t_0<t<\infty$ with $t_0=\frac{1}{\kappa}(1+\sqrt{1+\kappa^2})$,
and, thus, we have $\delta_0=b(t_0)$. It follows that, for $0\leq \delta\leq \delta_0$, $\delta=b(\eta)$ has
a unique solution $\eta$ with $ 0\leq \eta <t_0$.

By using the maximum principle for normal harmonic mappings, we prove
that a sequence of normal harmonic mappings $\{f_n\}$ converges to $0$ as
$n\rightarrow \infty$ in the unit disk under the condition that
$\max_{z\in C_n}|f_n(z)|$ converges to $0$, where $\{C_n\}$ is a
sequence of closed Jordan arcs with positive measure. Now, we state our next result
which is a generalization of \cite[Theorem 9.2]{Po1} for normal functions.

\bthm\label{thm4.2}
Suppose that $f_n$ are harmonic in $\mathbb{D}$ for $n\in \IN$, and
\be\label{eq4.9}
\sup_{z\in\mathbb{D}}(1-|z|^{2}){f_n}^{\#}(z)\leq\alpha<\infty, ~n\in\IN.
\ee
If there exist Jordan arcs $C_n\subset\mathbb{D}$ such that
\be\label{eq4.10}
{\rm diam}\,(C_n) =\sup_{z,w\in C_n}|z-w|\geq\gamma>0,~ n\in\IN,
\ee
and
\be\label{eq4.11}
\max_{z\in C_n}|f_n(z)|\to 0~\mbox{ as $n\to\infty$},
\ee
then $f_n(z)\to 0$ as $n\to\infty$, locally uniformly in $\mathbb{D}$.
\ethm

\bdefe\label{de1.3}
We say that a harmonic mapping $f$ in $\mathbb{D}$ has the {\it asymptotic value} $a\in\mathbb{C}$
at the point $\xi\in \mathbb{T}:=\{z:\, |z|=1\}$ if there exists a Jordan arc $\Gamma$ that
ends at $\xi$ and lies otherwise in $\mathbb{D}$ such that $f(z)\rightarrow a $ for $z\in \Gamma,\,z\rightarrow\xi$.

We call such an arc an {\it asymptotic path}. If $\Gamma=\{\xi
r:0\leq r\leq 1\}$, we call $a$ a {\it radial limit} (cf. {\cite{Po1}}).
\edefe

\bdefe\label{de1.4} A (symmetric) {\it Stolz angle} is a set of the
form
$$ A=\{z\in \mathbb{D}:
|\arg(1-\overline{\xi}z)|<(\pi/2)-\delta\} \;\; (0<\delta<\pi/2).
$$
That is, it is a sector with vertex $\xi$ and angle less than $\pi$ symmetric
to $[0, \xi]$. We say that $f$ has the {\it angular limit} $a$ at
$\xi\in \mathbb{T}$ if $f(z)\rightarrow a$ as $z\rightarrow\xi, \, z\in A$ and
for every Stolz angle $A$ at $\xi$ (cf. {\cite{Po1}}).
\edefe

By Definition \ref{de1.3}, an angular limit is a radial limit and
therefore is an asymptotic value. In the following theorem, we
show that the converse is true for normal harmonic mappings. Therefore
a normal harmonic mapping has at most one asymptotic value at any
given point $\xi \in \mathbb{D}$.

\bthm\label{thm4.3}
If the normal harmonic mapping $f$ has the asymptotic value $a$ at $\xi$,
then $f$ also has the angular limit $a$ at $\xi$.
\ethm

In Section \ref{sec-2}, we recall and also prove several lemmas which are useful to prove our main results.
In Section \ref{sec-3}, we present the proofs of the main theorems.

\section{Several Lemmas}\label{sec-2}

We begin this section with the following lemma which is a generalization of the corresponding one for analytic functions
due to  Marty (cf. \cite[p.~169]{Al}).

\blem\label{lem1} A class $\mathfrak{F}$ of harmonic mappings $f=h+\overline{g}$ in $\mathbb{D}$
is normal if $\{f^{\#}(z):\, f\in \mathfrak{F}\}$ (where $f^{\#}$ is defined in Definition \ref{de1.1})
is uniformly locally bounded.
\elem
\bpf Consider $\chi(f(z_1),f(z_2))$ defined as in \eqref{DPQ2-e q1} for $f(z_1)\neq \infty \neq f(z_2)$.
It is then easy to see that, followed by the stereographic projection, $f$ maps an arc $\gamma$ on an
image with length
\begin{eqnarray}
\nonumber L(\gamma)=\int_{\gamma}\frac{|df(z)|}{1+|f(z)|^{2}}
\leq \int_{\gamma}\frac{(|h'(z)|+|g'(z)|)\,|dz|}{1+|h(z)+\overline{g(z)}|^2}
=\int_{\gamma}{f^{\#}(z)\,|dz|}.
\end{eqnarray}
If $f^{\#}(z)\leq M$ on the segment between $z_1$ and $z_2$, where $M>0$ is  independent of $f$, then  we have
$$\chi(f(z_1),f(z_2))\leq \int_{\gamma}\rho(f)\,|dz|\leq M\int_{\gamma}|dz|=M|z_1-z_2|,
$$
which implies that harmonic mappings in $\mathfrak{F}$  are equicontinuous when $f^{\#}(z)$'s are locally bounded.
By Arzel\`{a}-Ascoli Theorem, the class $\mathfrak{F}$ is normal.
\epf

For a harmonic mapping $f=h+\overline{g}$ in $\ID$ such that $f(z_0)=0$ for some $z_0\in \mathbb{D}$, we have
the power series expansions of $h$ and $g$ in $|z-z_0|<1-|z_0|$ of the form
$$h(z)=a_0+\sum_{k=n}^{\infty}a_k(z-z_0)^k, \,\,\mbox{and}\,\,  g(z)=b_0+\sum_{k=m}^{\infty}b_k(z-z_0)^k,
$$
where $f(z_0)=a_0+\overline{b_0}=0$, $a_n\neq 0$ and $b_m\neq 0$. If $m\neq n$ or $m=n$ and $|a_n|\neq|b_m|$, we say that $f$ has a zero of order $\min \{m, n\}$ at $z_0$.

It is known that the zeros of a sense-preserving harmonic mapping are isolated (cf. \cite[p.~8]{Pe}).
We now recall the following lemma which is indeed the Hurwitz theorem for harmonic mappings.

\blem\label{lem2.1} {\rm(}\cite[p.~10]{Pe}{\rm )}
If $f$ and ${f_n}$ $(n\geq 1)$ are sense-preserving harmonic mappings in $\mathbb{D}$, and $\{f_n\}_{n\ge 1}$
converges locally uniformly to $f$, then $z_0\in\ID$ is a
zero of $f$ if and only if it is a cluster point of the zeros of the functions ${f_n}$ $(n\geq 1)$.
\elem

From Lemma \ref{lem2.1}, we observe that if $f$ has a zero of order $n$ at $z_0$
if and only if each small neighborhood of $z_0$ (small enough to contain no other zeros of $f$) contains
precisely $n$ zeros, counted according to multiplicity, of $f_n$ for every $n$ sufficiently large.
We say that $z=z_0$ is a multiple solution of $f(z)=\lambda$ if $z_0$ is a zero of order $n\geq 2$ of
$f(z)-\lambda$, that is $f(z_0)=h(z_0)+\overline{g(z_0)}=\lambda$, $|h'(z_0)|\neq0$ and $|g'(z_0)|\neq 0$.

Using Lemma \ref{lem2.1} and \cite[Corollary 3]{Hi}, we prove the following lemma.

\blem\label{lem2.2}
Let $f=h+\overline{g}$ be a sense-preserving harmonic mapping in $\mathbb{C}$ with $g(0)=0$. There are at
most four values of $\lambda$ for which all solutions of $f(z)=\lambda$ are multiple solutions.
\elem
\bpf
Let $f=h+\overline{g}$ be a sense-preserving harmonic mapping in $\mathbb{C}$ and  $\omega(z)=\frac{g'(z)}{h'(z)}$.
Then $|\omega(z)|<1$ in $\IC$ and thus, by Liouville's theorem,  $\omega(z)\equiv\alpha$ with the constant $|\alpha|<1$.
This gives
$$f(z)=h(z)+\overline{\alpha  h(z)-\alpha h(0)}.
$$
Now, for any number $\lambda$,  $f(z)=\lambda$ is equivalent to
$$h(z)=\frac{\lambda-\overline{\alpha\lambda}+  \overline{\alpha h(0)}-|\alpha|^2h(0)}{1-|\alpha|^2}.
$$
Thus if all solutions of $f(z)=\lambda$ are multiple solutions, then so do the last equation.
The converse is also true. By using \cite[Corollary 3]{Hi},
there are at most four values $\lambda^*$ for which all solutions of $h(z)=\lambda^*$ are multiple
solutions, which implies that there are at most four values $\lambda$ for which all solutions of
$f(z)=\lambda$ are multiple solutions.
\epf

\blem\label{lem3}{\rm (}\cite[Remark 1]{AHS}{\rm )} Let $\varphi$ be analytic in $\mathbb{D}$ and $|\varphi(z)|<1$.
If $f=h+\overline{g}$ is a  normal harmonic mapping in $\mathbb{D}$ and
$\sup_{z\in\mathbb{D}}(1-|z|^2)f^{\#}(z)=\alpha<\infty, $
then $F=f\circ\varphi=H+\overline{G}$ is also normal in $\mathbb{D}$, and
$\sup_{z\in\mathbb{D}}(1-|z|^2)F^{\#}(z)\leq\alpha
$
with equality if $\varphi\in \Aut (\mathbb{D})$.
\elem

Finally, we recall the identity theorem for harmonic mappings (\cite{AHS,Yo}).

\blem\label{lem4.3} Let $f$ be harmonic in a connected open set $D$. If $f(z)\equiv 0$ in some open subset $G\subset D$,
then $f(z)\equiv0$ in $D$.
\elem

\section{Proofs of theorems}\label{sec-3}

By using the method of proof of \cite[Theorem 1]{LP}, one can easily prove Theorem \ref{thm1}
but for the sake of completeness, we include the details.

\subsection{The proof of  Theorem \ref{thm1}}

 Suppose that $f$ is not normal. Then
there exists a sequence $\{z^{*}_n\}$ such that
\be\label{eq2.1}
(1-|z^{*}_n|^2)f^{\#}(z^{*}_n)\to \infty  ~\mbox{ as $n\to \infty$},
\ee
which also implies that $|z^{*}_n|\to 1 $ as $n\to \infty$.

Let $\{r_n\}$ be a sequence such that $|z^{*}_n|<r_n<1$ and
$$\left (1-\frac{|z^{*}_n|^{2}}{r^{2}_n}\right )f^{\#}(z^{*}_n)\to \infty  ~\mbox{ as $n\to \infty$}.
$$
Furthermore, we choose $\{z_n\}$ such that
\begin{eqnarray}
\nonumber M_n= \max_{|z|<r_n}\left (1-\frac{|z|^{2}}{r^{2}_n}\right)f^{\#}(z)
=\left (1-\frac{|z_n|^{2}}{r^{2}_n}\right)f^{\#}(z_n).
\end{eqnarray}
 Since $|z^{*}_n|<r_n$, it follows from \eqref{eq2.1} that $M_n\to \infty$ as $n\to \infty$.
If we set
\be\label{eq2.2}
\rho_{n}=\frac{1}{M_{n}}\left (1-\frac{|z_n|^{2}}{r^{2}_n}\right)=\frac{1}{f^{\#}(z_n)},
\ee
then we have
\begin{eqnarray}\label{eq2.3}
\frac{\rho_n}{1-|z_n|}\leq
\frac{\rho_n}{r_{n}-|z_n|}=\frac{r_n+|z_n|}{r^{2}_n M_n} \leq\frac{2}{r_n M_n}\to 0 ~\mbox{ as $n\to \infty$}.
\end{eqnarray}

Let $F_n(\zeta)=f(z_n+\rho_n \zeta),$ where $|\zeta|<R_n=\frac{1-|z_n|}{\rho_n}$.  From \eqref{eq2.3} we also note that
$R_n\to \infty$ as $n\to \infty$. It follows from \eqref{eq2.2}  that
\be\label{eq2.4}
F^{\#}_n(0)=\rho_n f^{\#}(z_n)=1.
\ee

We apply Lemma \ref{lem1} to show that the sequence $\{F_{n}(\zeta)\}$ is normal. If $|\zeta|\leq R\leq R_{n}$, then, by \eqref{eq2.2},
\begin{eqnarray}
\nonumber F^{\#}_n(\zeta)&=&\rho_n f^{\#}(z_n+\rho_n\zeta)\leq\frac{\rho_n M_n}{1-r^{-2}_n|z_n+\rho_n\zeta|^{2}}\\
\nonumber &\leq&\frac{r_n+|z_n|}{r_n+|z_n|-\rho_{n}R}\left (\frac{r_n-|z_n|}{r_n-|z_n|-\rho_n R}\right )
\end{eqnarray}
which, by \eqref{eq2.3}, tends to 1 as $n\to \infty$, for each fixed $R$.
Hence $\{F_{n}(\zeta)\}$ is a normal sequence.   We may assume that $\{F_{n}(\zeta)\}$ converges locally
uniformly in $\mathbb{C}$. Then, the limit function $F(\zeta)$ is harmonic in $\mathbb{C}$, and is non-constant because, by \eqref{eq2.4}, $F^{\#}(0)=1\neq 0$.

Next, we prove the necessary part of the theorem. Let $f$ be normal in $\mathbb{D}$. Again, we recall that the functions
$F_n(\zeta)$ given by
$$F_n(\zeta)=f(z_n+\rho_n \zeta)
$$
are defined for $|\zeta|<\frac{1-|z_n|}{\rho_n}$, and by \eqref{eq2.0}, we also have $\frac{\rho_n}{1-|z_n|}\to 0$
as $n\to \infty$, which implies  that $\frac{\rho_n}{1-|z_n|-\rho_{n}|\zeta|}\to 0$ as  $n\to \infty$,
for $|\zeta|<\frac{1-|z_n|}{\rho_n}$. Since
$$F_n^{\#}(\zeta)=\rho_{n}f^{\#}(z_{n}+\rho_{n}\zeta)
\leq\frac{\rho_n}{1-|z_n|-\rho_{n}|\zeta|}(1-|z_{n}+\rho_{n}\zeta|^{2})f^{\#}(z_{n}+\rho_{n}\zeta)
$$
and $f$ is normal, we have
$$ (1-|z_{n}+\rho_{n}\zeta|^{2})f^{\#}(z_{n}+\rho_{n}\zeta)<\infty.
$$
Therefore, $F_n^{\#}(\zeta)\to 0$ as $n\to \infty$ and thus,
$F^{\#}(\zeta)=0$ for all $\zeta\in \mathbb{C}$, so that $F(\zeta)$ is a constant. This completes the proof of Theorem \ref{thm1}.
$\hfill \Box$

\smallskip

\subsection{The proof of  Theorem \ref{thm2}}

Suppose that $f$ is a sense-preserving harmonic mapping in $\mathbb{D}$ which is not normal. By  Theorem \ref{thm1},
there exist sequences $\{z^*_n\}$ and $\{\rho_n\}$ with  $z^*_n\in \mathbb{D}$, $|z^*_n|\to 1$, $\rho_{n}>0$,
$\frac{\rho_n}{1-|z^*_n|}\to 0$ and a non-constant sense-preserving harmonic mapping $F$
in $\mathbb{C}$ such that the sequence $\{F_{n}\}$, $F_{n}(z)=f(z^*_{n}+\rho_{n}z)$, converges locally
uniformly to $F$ as $n\to \infty$.

Let $\lambda$ be any complex number, finite or infinite, for which the equation $F(t)=\lambda$ has a solution $z_0$
which is not a multiple solution, that is $F^{\#}(z_0)\neq0$. By Lemma \ref{lem2.1}, in each neighborhood of $z_0$ all but a
finite number of the functions $F_n$ assume the value $\lambda$. Thus there exists a sequence of points ${z_n}$ such that
$z_{n}\to z_{0}$  as $n\rightarrow \infty$, and $F_{n}(z_n)=\lambda$ for sufficiently large values of  $n$.  Also, since the convergence of  $\{F_n\}$
to $F$ is locally uniform, we have that $F^{\#}_{n}(z_n)\to F^{\#}(z_0)$. Letting $s_{n}=z^*_{n}+\rho_{n}z_{n}$, we get
that $F^{\#}_{n}(z_n)=\rho_{n}f^{\#}(s_n)$ so that
$$ f^{\#}(s_n)(1-|s_n|)
=F^{\#}_{n}(z_n)\frac{1-|z^*_n|}{\rho_n}\left (\frac{1-|s_n|}{1-|z^*_n|}\right ).
$$
Letting $n\to \infty$, we have that $F^{\#}_{n}(z_n)\to F^{\#}_{n}(z_0)$, $\frac{1-|z^*_n|}{\rho_n}\to \infty$, and
$\frac{1-|s_n|}{1-|z^*_n|}\to 1$ which imply that $f^{\#}(s_n)(1-|s_n|)\to \infty$ and hence, $(1-|s_n|^{2})f^{\#}(s_n)\to \infty$.

Now we have shown that if the equation $F(z)=\lambda$ has a solution which is not a multiple solution, then
$$ \sup_{z\in f^{-1}(\lambda)}(1-|z|^{2})f^{\#}(z)=\infty.
$$
However, by Lemma \ref{lem2.2}, there can be at most four values of $\lambda$ for which all solutions to the equation $F(z)=\lambda$
are multiple solutions. Thus, we have that if $f$ is a sense-preserving harmonic mapping in $\mathbb{D}$ such that $f$ is not normal,
then for each complex number $\lambda$, with at most four exceptions, we have
$$ \sup_{z\in f^{-1}(\lambda)}(1-|z|^{2})f^{\#}(z)=\infty.
$$
The proof is complete. $\hfill \Box$

\smallskip

\subsection{The proof of Theorem \ref{thm3}}

Because $f$ is normal, by assumption, we have that
$$ f^{\#}(z)\leq \frac{c_1(f)}{1-|z|^2}\leq \frac{c_1(f)}{1-|z|}.
$$
Let $\sigma=\chi(K, 2K)$, and let
$$ A=\min\left \{\frac{1}{2}, \frac{\sigma}{2c_1(f)}\right \}.
$$
Thus, if $z_0\in\mathbb{D}$ such that $|f(z_0)|\leq K$, then we have  $|z-z_0|\leq A (1-|z_0|)$ implies
$$ \chi(f(z),f(z_0))\leq \int_{L}f^{\#}(z)\,|dz|\leq \sigma,
$$
where $L$ is the line segment between $z$ and $z_0$. Also, $|z-z_0|<A(1-|z_0|)$ implies that $|f(z)|\leq 2K$
which in turn gives that $F$ is normal, where
$$F(z)=(f\circ \varphi)(z)=(h\circ \varphi)(z)+\overline{(g\circ \varphi)(z)}
=H(z)+\overline{G(z)}, \quad \varphi (z)=\frac{z_0-z}{1-\overline{z_0}z}.
$$
As $|z-z_0|<A(1-|z_0|)$ implies that $|f(z)|\leq 2K$, we have $|z|<A$ implies that $|F(z)|\leq 2K$.

Let $R={A}/{2}$. Then we have, for $z=re^{i\theta}$ with $r<R$,
$$F(z)=\frac{1}{2\pi}\int_{0}^{2\pi}\frac{R^2-r^2}{|Re^{it}-z|^2}F(Re^{it})\,dt.
$$
As
$$\frac{R^2-r^2}{|Re^{it}-z|^2}=\frac{Re^{it}}{Re^{it}-z} +\frac{\overline{z}}{Re^{-it}-\overline{z}} ~\mbox{ and }~
\frac{d^m}{dz^m} \left (\frac{a}{a-z}\right )=\frac{m!a}{(a-z)^{m+1}},
$$
we have
$$ H(z)=\frac{1}{2\pi}\int_{0}^{2\pi}\frac{Re^{it}}{Re^{it}-z}F(Re^{it})\,dt
$$
and thus,
\beq\label{eq2.6}
\nonumber|H^{(m)}(0)|&=&\frac{1}{2\pi}\left |\int_0^{2\pi}\frac{m!Re^{it}}{R^{m+1}e^{it(m+1)}}F(Re^{it})\,dt\right |
\leq  \frac{2 K m!}{ R^m}.
\eeq

\noindent
{\bf Claim 1.} $|h^{(m)}(z_0)|(1-|z_0|^2)^m<E'_m(f,K)$, where $E'_m(f,K)$ is a constant which depends only on $m$, $f$ and $K$.

Let us prove Claim 1 by the method of induction. As $H(z)=h(\varphi(z))$, we first consider
\be\label{DPQ2-eq2}
\varphi(z)-z_0= -(1-|z_0|^2)\frac{z}{1-\overline{z_0}z}=-(1-|z_0|^2)\sum_{k=1}^{\infty} (\overline{z_0})^{k-1}z^k
\ee
so that  $\varphi^{(n)}(0)= n!(|z_0|^2-1)(\overline{z_0})^{n-1}$, and compute that
$$H'(0) = h'(\varphi(0))\varphi '(0) =h'(z_0)(|z_0|^2-1)
$$
and
\beqq
H''(0)&=& h''(\varphi(0))(\varphi '(0))^2 + h'(\varphi(0))\varphi ''(0) \\
&=&h''(z_0)(|z_0|^2-1)^2 + 2h'(z_0)(|z_0|^2-1)\overline{z_0}.
\eeqq
In particular, by \eqref{eq2.6}, we have
$$|h'(z_0)|(1-|z_0|^2)<\frac{2 K }{ R}=E'_1(f,K)
$$
and
$$|h''(z_0)|(1-|z_0|^2)^2<2E'_1(f,K)+\frac{4K}{ R^2}=E'_2(f,K).
$$
Thus, the desired claim follows for $m=1,2$.

In order to apply the method of induction, we need to get an expression for $H^{(m)}(0)$ for $m\geq 3$ and for this,
we consider again \eqref{DPQ2-eq2} and
$$h(z)= \sum_{n=0}^{\infty}\frac{h^{(n)}(z_0)}{n!}(z-z_0)^n ~\mbox{ for $|z-z_0|<A(1-|z_0|)$}
$$
so that
$$H(z)= \sum_{n=1}^{\infty}\frac{H^{(n)}(0)}{n!}z^n ~\mbox{ for $\ds |z|<\frac{1+z_0}{1-|z_0|}A$.}
$$

 For integers $k$ and $n$ with $1\leq k\leq n$, let $B_1(n)=1$ and
$$B_k(n)=B_{k-1}(k-1)+B_{k-1}(k)+\cdots+B_{k-1}(n-1).
$$
It is easy to verify that
$$|B_k(n)|<(n-k+2)^{k-1} ~\mbox{ for $3\leq k\leq n$}.
$$
For $m=3$, we see that
\beqq
H'''(0)&=& h'''(\varphi(0))(\varphi '(0))^3 + 3h''(\varphi(0)) \varphi'(0)\varphi ''(0)+ h'(\varphi(0))\varphi'''(0) \\
&=&h'''(z_0)(|z_0|^2-1)^3+6h''(z_0)(|z_0|^2-1)^2\overline{z_0} +6h'(z_0)(|z_0|^2-1)(\overline{z_0})^{2}
\eeqq
and the claim for $m=3$ is easily seen to be true. Next, for $m\ge 4$, we have
\beq
\nonumber\frac{H^{(m)}(0)}{m!}&=&h'(z_0)(|z_0|^2-1)(\overline{z_0})^{m-1}+\frac{h^{''}(z_0)}{2!}(|z_0|^2-1)^2(\overline{z_0})^{m-2}(m-1)\\
&&
\nonumber +\sum_{k=3}^{m-1}\frac{h^{(k)}(z_0)}{k!}(|z_0|^2-1)^k(\overline{z_0})^{m-k}B_k(m)+\frac{h^{(m)}(z_0)}{m!}(|z_0|^2-1)^m.
\eeq

Now, we assume that the claim is true for $m=1,2, \ldots ,n-1$, and  show that it is also true for $m=n$. Indeed,
 using the last expression for $m=n$, we obtain that
\beqq
\left |\frac{h^{(n)}(z_0)}{n!}(|z_0|^2-1)^n\right | &<&\frac{2K}{R^n}+E'_1(f,K)+(n-1)E'_2(f,K)\\
&& \hspace {.5cm}+\sum_{k=3}^{n-1}B_k(n) E'_k(f,K) = E'_n(f,K).
\eeqq

By using the similar argument as that of Claim 1, we have
$$\left |\frac{g^{(m)}(z_0)}{m!}(|z_0|^2-1)^n\right | <E''_m(f,K),
$$
where $E''_m(f,K)$ is a constant which depends only on $m$, $f$ and $K$. Now, if we let $E_m(f,K)=E'_m(f,K)+E''_m(f,K)$, then
$$ (1-|z|^{2})^m(|h^{(m)}(z)|+|g^{(m)}(z)|)\leq E_{m}(f,K) ~\mbox{ for each $z\in \mathbb{D}$  }
$$
and such that $|f(z)|\leq K$. $\hfill \Box$

\smallskip

\subsection{The proof of  Theorem \ref{thm4}}
We first assume that $f$ is normal in $\mathbb{D}$. Assume the contrary of the assertion that there is a pair of sequences
$\{z_n\}$ and  $\{w_n\}$ of $\mathbb{D}$ such that $\rho(z_n, w_n)\to 0$ as $n\to \infty$, but
$\alpha=\lim_{n\to \infty}f(z_n)\neq\lim_{n\to \infty}f(w_n)=\beta$. Put
$$f_n=f\circ \phi_{z_n} ~\mbox{ with $\phi_{z_n}(z)=z_n+(1-|z_n|)z$ and $\ds u_n=\frac{w_n-z_n}{1-|z_n|}$}.
$$
As $f(z_n)=f_n(0)$ and $f(w_n)=f_n(u_n)$, we have $f_n(0)\to \alpha$ and $f_n(u_n)\to \beta$ as $n\to \infty$.
Since $\rho(z_n, w_n)\to 0$ as $n\to \infty$, we see that $\rho(z_n, w_n)< {1}/{2}$ for all sufficiently large $n$.
It follows that $|1-\overline{w_n}{z_n}|\leq 4(1-|z_n|)$ for all such $n$. Hence $|u_n|\leq 4\rho(z_n, w_n)\to 0$ as $n\to \infty$.
Furthermore,
\beq\label{eq2.7}
\nonumber\chi(f_n(0), f_n(u_n))&\leq& |u_n|\int_{0}^{1}f_n^{\#}(tu_n)\,dt=|u_n|\int_{0}^{1}f^{\#}(\phi_{z_n}(t u_n))|\phi_{z_n}'(t u_n)|\,dt\\
\nonumber&=&|w_n-z_n|\int_{0}^{1}f^{\#}(\phi_{z_n}(t u_n))\,dt\\
\nonumber&\leq& 2\left (\sup_{z\in \mathbb{D}} f^{\#}(z)(1-|z^2|)\right )\frac{|w_n-z_n|}{1-|z_n|} \to 0 ~\mbox{ as $n\to \infty$}
\eeq
which contradicts the fact that $\chi(\alpha, \beta)=0$.

Conversely, suppose that $\lim_{n\to \infty}f(z_n)=\lim_{n\to \infty}f(w_n)$ for all sequences
$\{z_n\}$ and  $\{w_n\}$ in $\mathbb{D}$ such that $\rho(z_n, w_n)\to 0$.

Let $\varphi_n\in \Aut (\mathbb{D})$,  $n=1, 2, \ldots$, and  $z_0\in \mathbb{D}$. Also, we assume that $\{w_n\}$ is a sequence of points in $\mathbb{D}$ such that $w_n\to z_0$ and $f(\varphi_n(w_n))$ converges to $\alpha$ for some $\alpha$.
For the sequence $\{z_n\}$ in $\mathbb{D}$, obviously, if $\rho(z_n, w_n)\to 0$ as $n\to \infty$, then $\rho(f(\varphi_n(z_n)), f(\varphi_n(w_n)))\to 0 $ as $n\to \infty$.
It follows that $f(\varphi_n(z_n))\to \alpha $ as $n\to \infty$, and then $f\circ \varphi_n$ is continuously convergent at $z_0$.
Since $z_0$ is an arbitrary point in $\mathbb{D}$,
we conclude that $f\circ \varphi_n$ is continuously convergent at each point of $\mathbb{D}$. Hence $\{f\circ \varphi_n\}$ is a normal family,
and thus, $\{h\circ \varphi_n\}$ and $\{g\circ \varphi_n\}$  are normal families. Therefore  $h$ and $g$ are normal (cf. \cite{No}).
By the assumption,  either $h$ or $g$ is bounded and thus, without loss of generality, we may assume that
$g$ is bounded, i.\,e., $|g(z)|\leq M$ in $\mathbb{D}$ for some $M>0$. For $z\in \mathbb{D}$ such that $|g(z)|<\frac{|h(z)|}{3}$, we have
\begin{eqnarray}
\nonumber(1-|z|^2)f^{\#}(z)&=&(1-|z|^2)\frac{|h'(z)|+  |g'(z)|}{1+\big |h(z)+\overline{g(z)}\big |^2}\\
\nonumber&\leq &(1-|z|^2)\frac{|h'(z)|}{1+\frac{1}{3}|h(z)|^2}+(1-|z|^2)|g'(z)|\\
\nonumber&\leq &3(1-|z|^2)h^{\#}(z)+(1+M)(1-|z|^2)g^{\#}(z)\\
\nonumber&<& \infty.
\end{eqnarray}
For $z\in \mathbb{D}$ such that $|g(z)|\geq \frac{|h(z)|}{3}$,
\begin{eqnarray}\nonumber (1-|z|^2)f^{\#}(z)&=&(1-|z|^2)\frac{|h'(z)|+|g'(z)|}{1+\big |h(z)+\overline{g(z)}\big|^2}\\
\nonumber&\leq &(1-|z|^2)|h'(z)|+(1-|z|^2)|g'(z)|\\
\nonumber&\leq &(1+9M^2)(1-|z|^2)h^{\#}(z)+(1+M)(1-|z|^2)g^{\#}(z)\\
\nonumber&<& \infty.
\end{eqnarray}
The preceding argument shows that $f$ is normal in $\mathbb{D}$. $\hfill \Box$

\smallskip
\begin{figure}
\begin{center}
\includegraphics[height=8cm, width=10.2cm, scale=1]{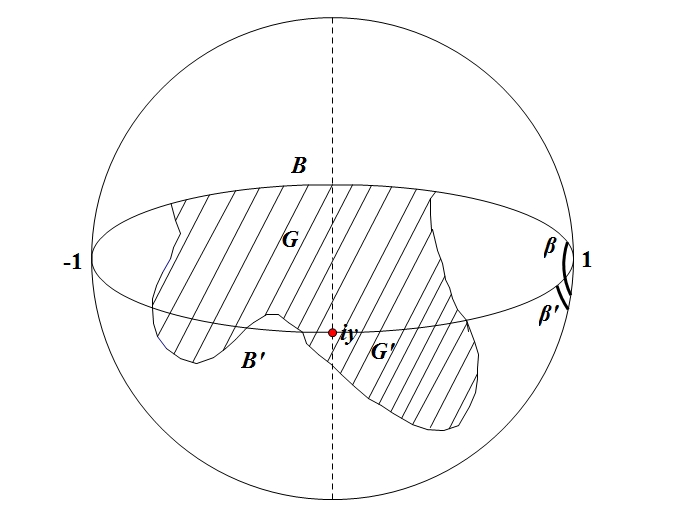}
\end{center}
\caption{Lens-shaped domain of angle
$\beta$ $(0<\beta<\pi)$ cut off from $\mathbb{D}$ by the circular arc $B$\label{maximum}  }
\end{figure}
\subsection{The proof of  Theorem \ref{thm4.1}}
By choosing a suitable  $\varphi\in \Aut (\mathbb{D})$
and replacing $f$ by $F=f\circ\varphi$, we assume that ${B}$ is a circular arc passing through $-1$ and $1$, and $G$ lies below the arc
${B}$. By Lemma \ref{lem3}, $F=f\circ\varphi=H+\overline{G}$ is also normal in $\mathbb{D}$, and
$$ \sup_{z\in\mathbb{D}}(1-|z|^2)F^{\#}(z)\leq\alpha .
$$
For $0<\beta'<\beta$, let $G'$ be the intersection of $G$ with the domain of angle $\beta'$ cut off
by the circular arc ${B'}$ through $\pm1$. See Figure \ref{maximum}.

Suppose that \eqref{eq4.2} does not hold. Since $|f(z)|\leq\delta<\eta$ for $z\in\partial G \backslash B$ by
\eqref{eq4.1}  and \eqref{eq4.3}, and since $G\subset \mathbb{D}$, there exists $\beta'$ such that $|f(z)|\leq\eta$ for
$z\in G'$. Let $\beta'~(0<\beta'<\beta)$ be the largest such number. Then
$$
\eta=\sup_{z\in G'}|f(z)|=|f(z_0)|
$$
for some $z_0\in B'\backslash\partial G$. By a further linear fractional transformation we may assume that $z_0=iy_0$, $y_0\in (-1,1)$, where
\be\label{eq4.5}
y_0=\tan \left (\frac{\beta'}{2}-\frac{\pi}{4}\right ).
\ee
Let
\be
\nonumber a(z)=|f(z)|\exp\left [\frac{b}{i}\log \Big (\frac{1+z}{1-z}\Big )+\frac{\pi b}{2}-2b\beta'\right ], \quad z\in G',
\ee
where $b=\frac{1}{\beta'}\log\frac{\eta}{\delta}>0$. It is known that every point in $\mathbb{D}$ lies
on one of the circular arcs that passes through $-1$, $iy$ and $1$ for some $y\in (-1,1)$, and on this circular arc
$$ \arg \left (\frac{1+z}{1-z}\right ) =\arg \left (\frac{1+iy}{1-iy}\right )= \arctan \left (\frac{2y}{1-y^2}\right )=
2\arctan y,
$$
from which it follows that $\exp[\frac{b}{i}\log\frac{1+z}{1-z}+\frac{\pi b}{2}-2b\beta']$
and $a(z)$ with $z\in\overline{G'}$ attains its maximum modulus on the boundary point $z_0=iy_0$. Since
$$ \max_{z\in B'\bigcap\partial G'}|a(z)|\leq\eta \exp(-b\beta')=\delta
$$
and $\partial G'\backslash B'\subset\partial G\backslash B$, we obtain from \eqref{eq4.1} that
$$ \sup_{z\in \partial G'\backslash B'}|a(z)|\leq \sup_{z\in \partial G'\backslash B'}|f(z)|\leq \sup_{z\in \partial G\backslash B}|f(z)|\leq\delta .
$$
Since $a(z)$ ($z\in\overline{G'}$) attains its maximum modulus on the boundary point $z_0=i y_0$,
it follows that $|a(z)|\leq\delta$ for $z\in G'$, so that
$$ \log|f(iy)|\leq \log\delta+2b\beta'-\frac{\pi b}{2}-2b\arctan y, ~iy\in G'.
$$
We have from \eqref{eq4.5} that $\beta'=\frac{\pi}{2}+2\arctan y_0$. Since $|f(iy_0)|=\eta$, we have
$$ \log|f(iy_0)|=\log\eta=\log\delta+b\beta'=\log\delta+2b\beta'-\frac{\pi b}{2}-2b\arctan y_0.
$$
Therefore
$$ \log|f(iy)|-\log|f(iy_0)|\leq-2b(\arctan y-\arctan y_0),~ iy\in G'.
$$
Letting $y\to y_0^{+}$ yields that
\be\label{eq4.6}
{\rm Re} \left (i\frac{h'(i y_0)-\overline{g'(i y_0)}}{f(i y_0)} \right )\leq\frac{-2b}{1+y_0^{2}}
=-\frac{2\log({\eta}/{\delta})}{\beta'(1+y_0^{2})}.
\ee
On the other hand, since $|f(iy_0)|=\eta$, it follows from \eqref{eq4.0} that
\be\label{eq4.7}
\left |\frac{h'(i y_0)-\overline{g'(i y_0)}}{f(i y_0)}\right |\leq\frac{|h'(i y_0)|+|g'(i y_0)|}{|f(i y_0)|}
\leq\frac{\alpha(1+\eta^{2})}{\eta(1-y_0^{2})}=\frac{\alpha(\eta+\eta^{-1})}{(1+y_0^{2})\sin\beta'}.
\ee
Hence \eqref{eq4.6}  and \eqref{eq4.7}  imply that
$$ \delta\geq\eta \exp \left [-\frac{k'}{2}\left (\eta+\frac{1}{\eta}\right )\right ] ~\mbox{ and }~ k'=\frac{\alpha\beta'}{\sin\beta'},
$$
which contradicts the fact that $k'<k$. The proof is complete. $\hfill \Box$

\smallskip

\subsection{The proof of  Theorem \ref{thm4.2}}
Suppose that the assertion is false. By \eqref{eq4.9} and Lemma \ref{lem1}, the sequence $\{f_n\}$ is normal in $\mathbb{D}$. Taking a subsequence we may therefore assume that
\be\label{eq4.12}
f_n(z)\to f(z)~\mbox{ as $n\to\infty$},
\ee
locally uniformly in $\mathbb{D}$, where $f$ is a harmonic mapping such that $f(z_0)\neq0$ for some $z_0\in\mathbb{D}$.

Now, we consider the first case that
\be\label{eq4.13}
\gamma_n=\inf\{|z|:z\in C_n\}\to 1~\mbox{ as $n\to\infty$}.
\ee
By \eqref{eq4.10}, there exist points $a_n,b_n\in C_n$ with $|a_n-b_n|= \gamma $. If $B_n$ denotes the circle through $a_n$ and
$b_n$ that is orthogonal to $\partial \ID$, then for sufficiently large values of $n$, $a_n$ and $b_n$ lie on different arcs of
$B_n^{*}=B_n\bigcap \{z_n:\, \gamma_n\leq|z_n|\leq1 \}$. Hence we can find a subarc $C'_n$ of $C_n$ that intersects each arc of $B_n^{*}$
exactly once. By \eqref{eq4.13} the subarc $B'_n$ of $B_n$ between the end points of $C'_n$ does not intersect $C'_n$ at any other point.
If $G_n$ is the inner domain of the Jordan curve $B'_n\bigcup C'_n$, then $\partial G_n=B'_n\bigcup C'_n\subset\mathbb{D}$,
which shows that $\overline{G_n}\subset\mathbb{D}$. Hence we obtain from \eqref{eq4.9}, \eqref{eq4.11} and Theorem \ref{thm4.1}
(with $\beta= {\pi}/{2}$) that
$$ \max_{z\in B'_n}|f_n(z)|\leq \max_{z\in\overline{G_n}}|f_n(z)|\to0~\mbox{ as $n\to\infty$}.
$$
Since $B'_n$ intersects the disk $\{z:\,|z|<r\}$ for some $r<1$ and for large values of $n$,
it therefore follows from \eqref{eq4.12} and Lemma \ref{lem4.3} that $f(z)\equiv0$, which is false.

In the case that \eqref{eq4.13} does not hold, $C_n$ intersects the closed disk $\{z:\,|z|\leq
r\}$ for some $r<1$ and for infinitely many values of $n$. Hence it follows from
\eqref{eq4.10}, \eqref{eq4.11}, \eqref{eq4.12} and Lemma \ref{lem4.3}
that $f(z)\equiv0$, which is again false.  The proof is complete. $\hfill \Box$

\smallskip

\subsection{The proof of  Theorem \ref{thm4.3}}

 Without loss of generality, we assume that $\xi=1$ and $a=0$. Suppose that $z_n\to1$ as $n\to\infty$ for $z_n\in A$,
where $A$ is a Stolz angle at $\xi$. First, we choose two real sequences $\{\xi_n\}$ and $\{y_n\}$, and $r<1$ such that
\be\label{eq4.15}
z_n=\varphi_n(i y_n), ~\varphi_n(s)=\frac{s+\xi_n}{1+\xi_n s},~ |y_n|\leq r, ~\mbox{ $\xi_n\to 1^{-}$ as $n\to\infty$}.
\ee
Obviously, $|z_n|<1$. The pre-image $\varphi_n^{-1}(\Gamma)$ of the asymptotic path $\Gamma$ intersects the
imaginary axis for sufficiently large values of $n$. Hence we can find a subarc $C_n$ of $\mathbb{D}\bigcap\varphi_n^{-1}(\Gamma)$
such that ${\rm diam}\,C_n\geq\frac{1}{2}$, ${\rm Re}\, z>0$ for $z\in C_n$, and there exists a sequence $\{w_n\}$ in the
arc $\varphi_n(C_n)\subset\Gamma$ such that $w_n\to1$ as $n\to\infty$. Then
\be\label{eq4.16}
\max_{s\in C_n}|f(\varphi_n(s))|=\max_{z\in\varphi_n(C_n)}|f(z)|\to 0~\mbox{ as $n\to\infty$}.
\ee
Since $f\circ\varphi_n$ is normal harmonic in $\mathbb{D}$ by Lemma \ref{lem1}, from \eqref{eq4.16}, we obtain that $f(\varphi_n(s))\to 0$
as $n\to \infty$, uniformly in $|s|\leq r$. Hence \eqref{eq4.15} shows that $f(z_n)\to 0 $ as $n\to\infty$.
The proof is complete.
$\hfill \Box$


\smallskip


\begin{thebibliography}{99}

\bibitem{Al} L. Ahlfors, Complex analysis,
\textit{MacGraw-Hill}, New York, 1953

\bibitem{AHS}H. Arbel\'{a}ez, R. Hern\'{a}ndez and W. Sierra, Normal harmonic mappings,
\textit{Monatsh Math} \textbf{190} (2019),  425–-439 



\bibitem{AMR}R. Aulaskari, S. Makhmutov and J. R\"{a}tty\"{a}, Results on meromorphic $\varphi-$normal functions,
\textit{Complex Var.  Elliptic Equ.}, {\bf 54} (2009), 855--863.

\bibitem{Co} F. Colonna,
The Bloch constant of bounded harmonic mappings,
\textit{Indiana Univ. Math. J.} {\bf 38}(1989), 829--840.

\bibitem{Pe}{ P. Duren}, Harmonic mappings in the plane,
\textit{Cambridge University Press}, New York, 2004.

\bibitem{Hi} E. Hille, Analytic function theory, V. II, \textit{AMS Chelsea publishing}, New York, 1962.

\bibitem{La0} P. Lappan, Some sequential properties of normal and non-normal functions with applications to
automorphic functions,
\textit{Comment. Math. Univ. St. Paul.}, {\bf 12} (1964), 41--57.

\bibitem{La1} P. Lappan, A criterion for a meromorphic function to be normal,
\textit{Comment. Math. Helv.}, {\bf 49} (1974), 492--495.

\bibitem{La2} P. Lappan, The spherical derivative and normal functions,
\textit{Ann. Acad. Sci. Fenn. Ser. A I Math.}, {\bf 3} (1977), 301--310.

\bibitem{LV} O. Lehto and K.~I. Virtanen, Boundary behaviour and normal meromorphic functions,
\textit{Acta Math.}, {\bf 97} (1957), 47--65.

\bibitem{Le} H. Lewy, On the non-vanishing of the Jacobian in certain one-to-one mappings,
\textit{Bull. Amer. Math. Soc.} \textbf{42}(1936), 689--692.

\bibitem{LP} A. J. Lohwater and Ch. Pommerenke, On normal meromorphic functions,
\textit{Ann. Acad. Sci. Fenn. Ser. A I Math.}, {\bf 550} (1973), 12 pp.

\bibitem{No} K. Noshiro, Contributions to the theory of meromorphic functions in the unit circle,
\textit{J. Fac. Svi. Hokkaido Univ.}, {\bf 7} (1938), 149--159.

\bibitem{Po} C. Pommerenke, On normal and automorphic functions,
\textit{Michigan Math. J.}, (1974), 193--202.

\bibitem{Po1}{Ch. Pommerenke}, Univalent functions.
\textit{Vandenhoeck $\&$ Ruprecht}, G\"{o}ttingen, 1975


\bibitem{Ya} S. Yamashita, On normal meromorphic functions,
\textit{Math. Z.}, {\bf 141} (1975), 139--145.

\bibitem{Yo} K.Yosida, On a class of meromorphic functions,
\textit{Proc. Phys. Math. Soc. Japan} 3. Ser. {\bf 16} (1934), 227--235.









\end{thebibliography}
\end{document}